\documentclass[12pt]{article}
\usepackage{url}
\usepackage{enumerate,amsfonts, graphicx, amsmath}
\usepackage{graphicx}
\usepackage{epsfig}
\usepackage{epstopdf}
\usepackage{multirow,color}

\title{A unified convergence bound for
conjugate gradient and accelerated gradient\thanks{Supported
in part by a grant from the U.~S.~Air Force Office of Scientific Research and in part
by a Discovery Grant from the Natural Sciences and Engineering Research Council
(NSERC) of Canada.}}
\author{Sahar Karimi\thanks{Department of Combinatorics \& Optimization,
University of Waterloo, 200 University Ave.~W., Waterloo, ON, N2L 3G1,
Canada, {\tt sahar.karimi@gmail.com}.} \and 
Stephen Vavasis\thanks{Department of Combinatorics \& Optimization,
University of Waterloo, 200 University Ave.~W., Waterloo, ON, N2L 3G1,
Canada, {\tt vavasis@uwaterloo.ca.}}}

\renewcommand\b{{\bf b}}
\renewcommand\d{{\bf d}}

\newcommand\p{{\bf p}}

\renewcommand\r{{\bf r}}
\newcommand\R{{\mathbb{R}}}
\newcommand\s{{\bf s}}
\renewcommand\t{{\bf t}}
\renewcommand\u{{\bf u}}

\newcommand\w{{\bf w}}
\newcommand\x{{\bf x}}
\newcommand\y{{\bf y}}

\newcommand\eps{\epsilon}

\newtheorem{theorem}{Theorem}
\newcommand{\eref}[1]{$(\ref{#1})$}

\begin{document}
\maketitle
\begin{abstract}
Nesterov's accelerated gradient
method for minimizing a smooth strongly convex function $f$ is
known to reduce $f(\x_k)-f(\x^*)$ by a factor 
of $\eps\in(0,1)$ after $k\ge O(\sqrt{L/\ell}\log(1/\eps))$ iterations, where
$\ell,L$ are the two parameters of smooth strong convexity.   Furthermore,
it is known that this is the best possible complexity in the function-gradient oracle
model of computation.  The method of linear conjugate gradients (CG) 
also satisfies the same
complexity bound in the special case of strongly convex quadratic functions,
but in this special case it is faster than the accelerated
gradient method.

Despite similarities in the algorithms and their
asymptotic convergence rates, the conventional analyses of the two methods are
nearly disjoint.  The purpose of this note is provide a single
quantity that decreases on every step at the correct rate for both
algorithms.
Our unified
bound is based on a potential similar to the potential
in Nesterov's original analysis.

As a side benefit of this analysis, we provide a direct proof that
conjugate gradient converges in $O(\sqrt{L/\ell}\log(1/\eps))$ iterations.
In contrast, the traditional indirect proof first establishes this
result for the Chebyshev algorithm, and then relies on optimality
of conjugate gradient to show that its iterates are at least as
good as Chebyshev iterates.  To the best of our knowledge, ours
is the first direct proof of the convergence rate of linear
conjugate gradient in the literature.
\end{abstract}


\section{Conjugate gradient}
The method of conjugate gradients  (CG)
was introduced by Hestenes and Stiefel \cite{hestenesstiefel}
for 
minimizing strongly convex quadratic functions
of the form $f(\x)=\x^TA\x/2-\b^T\x$, where $A$ is
a symmetric positive definite matrix.
We refer to this
algorithm as ``linear conjugate gradients.''

There is a significant body of work on
gradient methods for more general
smooth, strongly convex functions.  We say that a 
differentiable convex function $f:\R^n\rightarrow\R$ is
{\em smooth, strongly convex} 
\cite{hiriarturrutylemarechal}
if there exist two scalars $L\ge\ell>0$ such that
for all $\x,\y\in\R^n$,
\begin{equation}
\ell\Vert\x-\y\Vert^2/2\le f(\y)-f(\x)-\nabla f(\x)^T(\y-\x)\le L\Vert\x-\y\Vert^2/2.
\label{eq:strconvdef}
\end{equation}
This is equivalent to assuming convexity and lower and upper
Lipschitz constants on the gradient:
$$\ell\Vert\x-\y\Vert \le \Vert\nabla f(\x)-\nabla f(\y)\Vert \le
L\Vert\x-\y\Vert.$$

Nemirovsky and Yudin \cite{NemYud83}
proposed
a method for minimizing smooth strongly convex
functions requiring
$k=O(\sqrt{L/l}\log(1/\eps))$ iterations to produce an iterate 
$\x_k$ such that $f(\x_k)-f(\x^*)\le\eps(f(\x_0)-f(\x^*))$, where
$\x^*$ is the optimizer (necessarily unique under the assumptions made).  
A drawback of their method is that it requires an expensive two-dimensional
optimization on each iteration.  Nesterov \cite{Nesterov:k2} proposed
another method, nowadays
known as the ``accelerated gradient method,'' 
which achieves the same optimal complexity that requires a single
function and gradient evaluation on each iteration.

The accelerated gradient method, although optimal in theory, can be slow in practice.
For example, in the case of quadratic function, computational testing shows
that it is substantially slower  than linear conjugate
gradients.
In the special case of strongly convex quadratic functions, 
the conjugate
gradient has already been known to satisfy the same asymptotic
bound since the 1960s.

Although the two methods satisfy the same asymptotic bound,
the analyses of the two methods are completely
different.  In the case of accelerated gradient, there are
two analyses by Nesterov \cite{Nesterov:k2, Nesterov:book}.
There is also a recent analysis of a variant of accelerated
gradient \cite{bubeck}, which views it as a kind of ellipsoid method.
(This variant apparently requires exact line search.)  

In the case of linear
conjugate gradient, we are aware of no direct analysis of the algorithm.
By ``direct,'' we mean an analysis of $f(\x_{k})-f(\x^*)$ using the
recurrence inherent in CG.  Instead, the standard
analysis proves that another iterative method, for
example Chebyshev iteration \cite{GVL} or the heavy-ball iteration
\cite{Polyak,Bertsekas} achieves reduction of 
$\left(1-O(\sqrt{\ell/L})\right)$ per iteration.  Then one appeals
to the optimality of the CG iterate in the Krylov space generated
by all of these methods to claim that the CG iterate must be at
least as good as the others.

This paper is devoted to establishing
a one-step convergence bound that applies to both methods.
The one-step convergence bound has the form $\Psi_{k+1}\le \Psi_k/(1+\sqrt{\ell/L})$
for $k=1,2,\ldots$,
where $\Psi_k$ is a potential 
defined by \eref{eq:potdef}. This potential involves
both certain perturbed distance from the current iterate to the optimizer and
the objective function residual.
It should be noted that
for the accelerated gradient method, neither the
sequence $\Vert \x_k-\x^*\Vert$ nor $f(\x_k)-f(\x^*)$ is monotonically
decreasing with $k$.  Both of these sequences decrease monotonically for 
conjugate gradient (refer to
\eref{eq:cgdesc} and
\eref{eq:hs6:5} below), but neither decreases at the rate
$1/(1+\sqrt{\ell/L})$ on every step.  Instead, the rate of
decrease of these quantities (both in theory and in practice) 
is erratic.  Thus, it is not obvious that there is
a combination of these quantities that decreases at the proper rate
on every iteration for both algorithms.

The $k=0$ case of \eref{eq:potdef} is
$$(\ell/2)\Psi_0=(\ell/2)\Vert\x_0-\x^*\Vert^2+(f(\x_0)-f(\x^*)),$$
while $(\ell/2)\Psi_k\le f(\x_k)-f(\x^*)$.  Furthermore,
we show below that $\Psi_{k+1}\le \Psi_k/(1+\sqrt{\ell/L})$ 
for $k=1,2,\ldots$ and $\Psi_1\le \Psi_0$.
The consequence of all these bounds is the following theorem.
\begin{theorem}
Let $f(\x)$ be a strongly convex smooth function with
convexity parameters $\ell,L$.  Then the accelerated
gradient method
produces a sequence of iterates $\x_k$
such that
\begin{equation}
f(\x_k)-f(\x^*)\le C_0\left(1+\sqrt{\frac{\ell}{L}}
\right)^{-(k-1)},
\label{eq:fxkbd}
\end{equation}
where $\x^*$ is
the (necessarily unique) optimizer and
$C_0=(\ell/2)\Vert\x_0-\x^*\Vert^2+f(\x_0)-f(\x^*)$.
When applied to a quadratic function, the conjugate
gradient method produces a sequence satisfying this bound.
\label{mainthm}
\end{theorem}
Note that $(1+\sqrt{\ell/L})^{-1}\le (1-\sqrt{\ell/(4L)})$, so 
\eref{eq:fxkbd} implies the usual theorem except for a constant
factor.
Note that in the $k=0$ case, we establish only $\Psi_1\le\Psi_0$ instead
of the stronger $\Psi_{k+1}\le\Psi_k/(1+\sqrt{\ell/L})$, which is valid
for $k=1,2,\ldots$.
This explains why 
the exponent in Theorem~\ref{mainthm}
is $k-1$ rather than $k$.

In Section~\ref{sec:nest}, we review the accelerated gradient
method.  In Section~\ref{sec:cg}, we
review conjugate gradient and convergence bound.
In Section~\ref{sec:unifalg} we propose a single algorithmic
framework that unifies both algorithms.  Then, in the
main technical sections of this article, 
Sections~\ref{sec:unifana}--\ref{sec:cvgcg},
we present our unified analysis of the two algorithms, which is
an extension of the
potential-function approach used in Nesterov's original analysis.
Finally, in Section~\ref{sec:discussion}, we discuss some possible
consequences and future directions made possible by the
unified bound.

\section{Accelerated gradient method}
\label{sec:nest}

Following the treatment in his book \cite{Nesterov:book}
Nesterov's accelerated method can be described as follows. 
Given a strongly convex $f:\R^n\rightarrow\R$ with strong
convexity parameters $L,\ell$, one uses the recurrence:
\begin{align}
& \mbox{\bf Accelerated gradient method} \notag\\
& \x_0:=\mbox{arbitrary} \notag \\
&\mbox{for }k:=0,1,2,\ldots \notag  \\
&\hphantom{\mbox{for }} \y_{k+1}:=\x_k+\theta_k\s_k\label{eq:yupd}\\
&\hphantom{\mbox{for }} \x_{k+1}:=\y_{k+1}-\nabla f(\y_{k+1})/L  \\
&\hphantom{\mbox{for }} \s_{k+1}:=\x_{k+1}-\x_k \\
& \mbox{end} \notag
\end{align}
In \eref{eq:yupd} when $k=0$, $\s_0$ is undefined and hence we
define $\y_1:=\x_0$, and thus $\theta_0=0$.  For $k\ge 1$,
several choices of $\theta_k$ are valid; our analysis
uses 
\begin{equation}
\theta_k=\frac{\sqrt{L}-\sqrt{\ell}}{\sqrt{L}+\sqrt{\ell}}.
\label{eq:agthetak}
\end{equation}

\section{Conjugate gradient method}
\label{sec:cg}

The conjugate gradient method for 
minimizing $f(\x)=\x^TA\x/2 -\b^T\x$, 
where
$A$ is a symmetric positive definite matrix, 
is due to Hestenes and Stiefel
\cite{hestenesstiefel}
and is as follows.
\begin{align}
& \mbox{\bf LCG method} \notag \\
& \x_0:=\mbox{arbitrary} \notag \\
& \r_0:=\b-A\x_0 \notag \\
& \mbox{for } k:=0,1,2,\ldots, \notag \\
&\displaystyle \hphantom{\mbox{for }} 
\beta_{k+1}:= \frac{\r_{k}^T\r_{k}}{\r_{k-1}^T\r_{k-1}} \label{eq:betadef}\\
&\hphantom{\mbox{for }} 
\p_{k+1} :=  \beta_{k+1}\p_{k}+\r_{k} \label{eq:pupd} \\
&\displaystyle \hphantom{\mbox{for }} 
\alpha_{k+1}:= \frac{\r_{k}^T\r_{k}}{\p_{k+1}^TA\p_{k+1}} \label{eq:alphadef} \\
&\hphantom{\mbox{for }} 
\x_{k+1} := \x_{k}+\alpha_{k+1}\p_{k+1} \label{eq:xupd} \\
&\hphantom{\mbox{for }} 
\r_{k+1} := \r_{k}-\alpha_{k+1} A\p_{k+1} \label{eq:rupd} \\
&\mbox{end} \notag
\end{align}

When $k=0$, $\r_{k-1}$ is undefined.  Hence we disregard \eref{eq:betadef}
for specifying $\beta_1$ and instead take $\beta_1=0$, which implies
$\p_1=\r_0$ in \eref{eq:pupd}.  It is
apparent from this recurrence that
$\r_k=\b-A\x_k=-\nabla f(\x_k)$ for all $k$.
Here are two other well known relationships 
from \cite{hestenesstiefel}:
\begin{align}
\p_k^T\r_k& =0, &&\mbox{(HS 5:3c)}\label{eq:prorth} \\
\frac{1}{\alpha_k} &\in 
[\lambda_{\min}(A),\lambda_{\max}(A)]. && \mbox{(HS 5:12)}
\label{eq:alphaRQ}
\end{align}
Several monographs explain the method in detail from different
points of view
including Golub and Van Loan \cite{GVL}, Trefethen and Bau \cite{TB},
Greenbaum \cite{Greenbaum} and Liesen and Strakos \cite{LiesenStrakos}.


The best-known theorem regarding the convergence rate of conjugate 
gradient is due to Daniel \cite{Daniel} (but see
\cite{LiesenStrakos} for a more comprehensive perspective):
\begin{theorem}
For the above iteration,
$$f(\x_k)-f(\x^*)\le
4\left(\frac{1-\sqrt{\ell/L}}{1+\sqrt{\ell/L}}\right)^{2k}(f(\x_0)-f(\x^*)).$$
\end{theorem}

Daniel's proof and all others known to us use the following line
of reasoning.  First, Daniel uses a known result that the Chebyshev
method satisfies the bound above.  Then he relies on the fact that
the Chebyshev iterate $\x_k^{\rm Ch}$ lies in the affine space
$\x_0+{\rm span}\{\b, A\b,\ldots,A^{k-1}\b\}$.  On the other hand,
conjugate gradient is known to produce the vector $\x_k$ that
is the optimizer of $f$ over this affine space.  
Therefore, the
conjugate gradient iteration produces at least the same amount
of reduction in $f$. 
The analysis of the convergence rate
of conjugate gradient developed below
does not rely on the optimality with respect to the Krylov space.

Daniel's theorem is tight in the sense that for every choice of
$0<\ell<L$ and $k$, there a matrix $A$ and starting vector $\b$ such
that the bound in the theorem is achieved to within constant factors.
This follows from a much more general result of Nesterov \cite{Nesterov:book},
which states that the bound in Daniel's theorem is the best possible bound for
any algorithm that uses the function-gradient oracle model.  Linear conjugate
gradient applied to convex quadratic functions is a member of this
class of algorithms.
However, for particular choices of $A$, much better behavior may
be observed from linear conjugate gradient.

\section{Unified algorithm}
\label{sec:unifalg}

In this section, we consider the following iterative framework, which 
has three sequences of scalar parameters, $\theta_k$, $\nu_k$ and $\pi_k$
for $k=0,1,\ldots$.

\begin{align}
& \mbox{\bf Unified framework} \notag\\
& \x_0:=\mbox{arbitrary} \notag \\
&\mbox{for }k:=0,1,2,\ldots \notag  \\
&\hphantom{\mbox{for }} \y_{k+1}:=\x_k+\theta_k\s_k\label{eq:yupd2}\\
&\hphantom{\mbox{for }} \x_{k+1}:=\x_k+\nu_k\s_k-\pi_k\nabla f(\y_{k+1}) \label{eq:xupd2} \\
&\hphantom{\mbox{for }} \s_{k+1}:=\x_{k+1}-\x_k \label{eq:sformula}\\
& \mbox{end} \notag
\end{align}
When $k=0$, we leave $\s_0$ undefined and
take $\y_1:=\x_0$ in \eref{eq:yupd2} and $\x_1:=\x_0-\pi_0\nabla f(\y_1)$
in \eref{eq:xupd2}.  This in turn means that we start with $\nu_0=\theta_0=0$.

It is straightforward to observe that the accelerated gradient method is
a special case of the unified framework if we make the identification
$$\nu_k^{\rm AG}\equiv \theta_k^{\rm AG}\equiv 
\frac{\sqrt{L}-\sqrt{\ell}}{\sqrt{L}+\sqrt{\ell}}$$
for $k\ge 1$
and $\pi_k^{\rm AG}=1/L$ for all $k\ge 0$.

The LCG method can be derived as a special case of the
unified framework as follows.  First, take
$\theta_k^{\rm CG}\equiv 0$ so that $\y_{k+1}\equiv\x_k$ for all $k$.
Comparing
\eref{eq:sformula} and \eref{eq:xupd} we see that
\begin{equation}
\s_k=\alpha_k\p_k. \label{eq:sprel}
\end{equation}
Substituting \eref{eq:pupd} into
\eref{eq:xupd} yields
\begin{align*}
\x_{k+1}&=\x_k+\alpha_{k+1}(\beta_{k+1}\p_k+\r_k) \\
&=\x_k+\alpha_{k+1}\left(\frac{\beta_{k+1}}{\alpha_k}\s_k-\nabla f(\x_k)\right).
\end{align*}
We recover this recurrence if we take
\begin{align}
\nu_k^{\rm CG} &\equiv \frac{\alpha_{k+1}\beta_{k+1}}{\alpha_k}, 
&&k=1,2,\ldots,\label{eq:nucg}\\
\pi_k^{\rm CG} &\equiv \alpha_{k+1}, &&k=0,1,2\ldots \label{eq:picg}
\end{align}
in \eref{eq:xupd2}.

In both LCG and accelerated gradient, the parameters satisfy the following
relationships, which we assume for the rest of this paper:
\begin{equation}
\nu_k\ge \theta_k\ge 0\>(k=0,1,\ldots);\quad \nu_k>0\>(k=1,2,\ldots);\quad 
\pi_k> 0\>(k=0,1,\ldots).
\label{eq:paramsign}
\end{equation}

\section{A potential for both algorithms}
\label{sec:unifana}

In this section we propose the common potential for
both algorithms that decreases on every iteration.
The main result we establish is:
\begin{equation}
C\Psi_{k+1} \le \Psi_k
\label{eq:potdecr}
\end{equation}
where 
$$C = 1+\sqrt{\ell/L}$$
and
$\Psi_k$ is
a potential at step $k$:
\begin{equation}
\Psi_k = \Vert\w_k\Vert^2 + \frac{2}{\ell} (f(\x_k)-f(\x^*))
\label{eq:potdef}
\end{equation}
Here,
$$\w_k= \x_k+\rho_k \s_k-\x^*,$$
where
$\rho_k$, $k=0,1,\ldots$, is an additional sequences of scalars defined below 
(see \eref{eq:rhodefag} and \eref{eq:rhodeflcg}),
 $\ell$ is the lower strong-convexity
parameter, $L$ is the upper parameter, and $\x^*$ is the minimizer of $f$.
In fact, in the case of accelerated gradient, a slightly
stronger bound of 
$$C = 1 + \frac{1}{\sqrt{L/\ell}-1}$$
is established.
In the case $k=0$, we define $\w_0=\x_0-\x^*$ (hence $\rho_0=0$).

A potential involving
these two terms was proposed in \cite{Nesterov:k2}, and our analysis
may therefore be regarded as a variant of Nesterov's
technique.  (In \cite{Nesterov:k2}, only the second term 
of $\Psi_k$ is updated 
by a scalar from one iteration to the next.)

\section{Analysis of accelerated gradient}

We start by rewriting $\w_{k+1}$ and $\w_k$ in terms of 
$\y_{k+1}$, $\s_k$ and $\x^*$:
\begin{align*}
\w_k &= \x_k+\rho_k\s_k-\x^*  \\
&=\y_{k+1}+(\rho_k-\theta_k)\s_k -\x^*, \quad \mbox{(by \eref{eq:yupd2})}
\end{align*}
and
\begin{align*}
\w_{k+1} & = \x_{k+1}+\rho_{k+1} \s_{k+1}-\x^*   \\
& = (1+\rho_{k+1})\x_{k+1}-\rho_{k+1}\x_k-\x^* \quad \mbox{(by \eref{eq:sformula})} \\
& = (1+\rho_{k+1})(\y_{k+1}+(\nu_k-\theta_k)\s_k-\pi_k\nabla f(\y_{k+1}))
-\rho_{k+1}\x_k-\x^* \quad \mbox{(by \eref{eq:yupd2} and \eref{eq:xupd2})} \\
&= \y_{k+1}+\rho_{k+1}(\y_{k+1}-\x_k)+(1+\rho_{k+1})((\nu_k-\theta_k)\s_k-\pi_k\nabla f(\y_{k+1}))-\x^* \\
& = \y_{k+1} +(\rho_{k+1}\theta_k+(1+\rho_{k+1})(\nu_k-\theta_k))\s_k  -(1+\rho_{k+1})\pi_k\nabla f(\y_{k+1})-\x^*
\quad \mbox{(by \eref{eq:yupd2})} \\
& = 
\y_{k+1}-\x^* +((1+\rho_{k+1})\nu_k-\theta_k)\s_k  -(1+\rho_{k+1})\pi_k\nabla f(\y_{k+1}).  
\end{align*}

We now let $\xi=\sqrt{C}$, which implies that
the first term of $C\Psi_{k+1}-\Psi_k$ is of the  form:
\begin{equation}
\Vert\xi\w_{k+1}\Vert^2 - \Vert\w_{k}\Vert^2
= (\xi\w_{k+1}-\w_k)^T(\xi\w_{k+1}+\w_k).
\label{eq:potdif1}
\end{equation}
We expand the two factors separately using the previously
developed expressions for $\w_{k+1}$ and $\w_k$:
\begin{align*}
\xi\w_{k+1}-\w_k & = 
\xi(\y_{k+1}-\x^* +((1+\rho_{k+1})\nu_k-\theta_k)\s_k  -
(1+\rho_{k+1})\pi_k\nabla f(\y_{k+1})) \\
& \hphantom{=}\quad\mbox{}
- 
(\y_{k+1}+(\rho_k-\theta_k)\s_k -\x^*) \\
&\equiv \t_1+\t_2-\t_3
\end{align*}
where
\begin{align*}
\t_1 &= (\xi-1)(\y_{k+1}-\x^*), \\
\t_2 &= (\xi((1+\rho_{k+1})\nu_k-\theta_k)-(\rho_k-\theta_k))\s_k, \\
\t_3 & =  \xi(1+\rho_{k+1})\pi_k\nabla f(\y_{k+1}).
\end{align*}

Here, the vectors $\t_1,\t_2,\t_3$ also depend on iteration $k$, but we omit
writing this dependence since $k$ is fixed for this part of the analysis.
Similarly,
\begin{align*}
\xi\w_{k+1}+\w_k & = 
\xi(\y_{k+1}-\x^* +((1+\rho_{k+1})\nu_k-\theta_k)\s_k  -(1+\rho_{k+1})\pi_k\nabla f(\y_{k+1})) \\
& \hphantom{=}\quad\mbox{}
+ 
(\y_{k+1}+(\rho_k-\theta_k)\s_k -\x^*) \\
&\equiv \u_1+\u_2-\u_3
\end{align*}
where
\begin{align*}
\u_1 &= (\xi+1)(\y_{k+1}-\x^*), \\
\u_2 &= (\xi((1+\rho_{k+1})\nu_k-\theta_k)+(\rho_k-\theta_k))\s_k, \\
\u_3 &=  \xi(1+\rho_{k+1})\pi_k\nabla f(\y_{k+1}).
\end{align*}
Thus, \eref{eq:potdif1}
is rewritten $(\t_1+\t_2-\t_3)^T(\u_1+\u_2-\u_3)$. This expansion  contains nine
terms.
Writing these and gathering like terms (and noting the simple identity
$(a-b)(c+d)+(a+b)(c-d)=2ac-2bd$) yields:
\begin{align}
\t_1^T\u_1 & = (\xi^2-1)\Vert\y_{k+1}-\x^*\Vert^2, \label{eq:tu1}\\
\t_1^T\u_2+\t_2^T\u_1 &= 2(\xi^2((1+\rho_{k+1})\nu_k-\theta_k)
-(\rho_k-\theta_k))(\y_{k+1}-\x^*)^T\s_k, \label{eq:tu2}\\
\t_2^T\u_2 &=
(\xi^2((1+\rho_{k+1})\nu_k-\theta_k)^2-(\rho_k-\theta_k)^2)
\Vert\s_k\Vert^2, \label{eq:tu3}\\
-\t_1^T\u_3-\t_3^T\u_1 &=2\xi^2(1+\rho_{k+1})\pi_k(\x^*-\y_{k+1})^T
\nabla f(\y_{k+1}), \label{eq:tu4}\\
-\t_2^T\u_3-\t_3^T\u_2 &=-2\xi^2((1+\rho_{k+1})\nu_k-\theta_k)(1+\rho_{k+1})
\pi_k \s_k^T\nabla f(\y_{k+1}), \label{eq:tu5}\\
\t_3^T\u_3 &= \xi^2(1+\rho_{k+1})^2\pi_k^2\Vert\nabla f(\y_{k+1})\Vert^2.
\label{eq:tu6}
\end{align}

For accelerated gradients, we
use a constant value for
$\rho_k$ (independent of $k$) that is analogous to the
choice in \cite{Nesterov:k2}, namely,
\begin{equation}
\rho_k=\sqrt{L/\ell} - 1\quad\mbox{for $k=1,2,\ldots$}.
\label{eq:rhodefag}
\end{equation}
Assume for now that $k\ge 1$; 
the $k=0$ case is considered separately below.
The inner product $(\y_{k+1}-\x^*)^T\s_k$ in \eref{eq:tu2}
appears difficult to bound in the case of accelerated gradient, so
we define the scalar $\xi^2$ ($=C$) to ensure that the
term $\t_1^T\u_2+\t_2^T\u_1$ is zero, namely,
\begin{align}
\xi^2 &= \frac{\rho_k-\theta_k}
{(1+\rho_{k+1})\nu_k-\theta_k} \notag\\
&=
\frac{\sqrt{L/\ell}-1-(\sqrt{L/\ell}-1)/(\sqrt{L/\ell}+1)}
{(\sqrt{L/\ell}-1)\cdot(\sqrt{L/\ell}-1)/(\sqrt{L/\ell}+1)} \notag\\
&=
1+\frac{1}{\sqrt{L/\ell}-1}. \label{eq:xirecurag}
\end{align}
Note that this implies $C\ge 1 + \sqrt{\ell/L}$, so that 
\eref{eq:fxkbd} will be established for this choice of $\xi$.

Next, we rewrite the remaining terms of \eref{eq:tu1}--\eref{eq:tu6}
based on these choices for the scalars:
\begin{align}
\t_1^T\u_1 & = \frac{1}{\sqrt{L/\ell}-1}\Vert\y_{k+1}-\x^*\Vert^2, \label{eq:tu1ag}\\
\t_2^T\u_2 &=
-\frac{\sqrt{L/\ell}(\sqrt{L/\ell}-1)^2}{(\sqrt{L/\ell}+1)^2}
\Vert\s_k\Vert^2, \label{eq:tu3ag}\\
-\t_1^T\u_3-\t_3^T\u_1 &=\frac{2}{(\sqrt{L/\ell}-1)\ell}(\x^*-\y_{k+1})^T
\nabla f(\y_{k+1}), \label{eq:tu4ag}\\
-\t_2^T\u_3-\t_3^T\u_2 &=
-\frac{2}{\ell}\cdot\frac{\sqrt{L/\ell}-1}{\sqrt{L/\ell}+1}
\cdot
\s_k^T\nabla f(\y_{k+1}), \label{eq:tu5ag}\\
\t_3^T\u_3 &= \frac{1}{L^{1/2}\ell^{3/2}(\sqrt{L/\ell}-1)}\Vert\nabla f(\y_{k+1})\Vert^2.
\label{eq:tu6ag}
\end{align}

We analyze the sum of \eref{eq:tu4ag}, \eref{eq:tu5ag}, and \eref{eq:tu6ag} together:
\begin{equation}
-\t_1^T\u_3-\t_3^T\u_1
-\t_2^T\u_3-\t_3^T\u_2+
\t_3^T\u_3
=
\frac{2}{(\sqrt{L/\ell}-1)\ell}\cdot t_4
\label{eq:threeterm}
\end{equation}
where
\begin{equation}
t_4=(\x^*-\y_{k+1})^T
\nabla f(\y_{k+1})
-\frac{(\sqrt{L/\ell}-1)^2}{\sqrt{L/\ell}+1}
\s_k^T\nabla f(\y_{k+1})
+\frac{1}{2\sqrt{L\ell}}
\Vert\nabla f(\y_{k+1})\Vert^2.
\label{eq:t4def}
\end{equation}

To analyze $t_4$ requires two more bounds.
First, by \eref{eq:strconvdef},
for any $\x\in\R^n$,
\begin{equation}
f(\x)\ge f(\y_{k+1})+\nabla f(\y_{k+1})^T(\x-\y_{k+1}) +\frac{\ell}{2}
\Vert\x-\y_{k+1}\Vert^2.
\label{eq:strcvx}
\end{equation}
We also need
a bound on the descent made per step.  
We use the well known bound
\begin{equation}
f\left(\y_{k+1})-f(\y_{k+1}-\nabla f(\y_{k+1})/L\right)\ge \Vert\nabla  f(\y_{k+1})
\Vert^2/(2L).
\label{eq:cvxdesc}
\end{equation}
This follows by writing the left-hand side 
$f(\y_{k+1})-f(\y_{k+1}-\d)$ 
as the line integral
$\int_0^1\nabla f(\y_{k+1}-t\d)^T\d\,dt$ for the particular choice
$\d=\nabla f(\y_{k+1})/L$, pulling out an additive
term of $\Vert\nabla f(\y_{k+1})\Vert^2/L$ from the integrand, 
and then applying
the Lipschitz condition.

Then the claimed bound is:
\begin{equation}
t_4\le 
(\sqrt{L/\ell}-1)(f(\x_k)-f(\x^*))-\sqrt{L/\ell}(f(\x_{k+1})-f(\x^*))
-\frac{\ell}{2}\Vert\y_{k+1}-\x^*\Vert^2.
\label{eq:t4bound}
\end{equation}
The following chain of inequalities starting
from \eref{eq:t4def} establishes \eref{eq:t4bound}:
\begin{align*}
t_4 
&= 
(\x^*-\y_{k+1})^T\nabla f(\y_{k+1}) 
+ (\sqrt{L/\ell}-1)(\x_k-\y_{k+1})^T\nabla f(\y_{k+1}) 
+ 
\frac{\Vert\nabla f(\y_{k+1})\Vert^2}{2\sqrt{L\ell}}
\\
&\qquad\qquad
\mbox{(by \eref{eq:yupd2})} \\
&\le
f(\x^*)-f(\y_{k+1})-\frac{\ell}{2}\Vert \y_{k+1}-\x^*\Vert^2
+ (\sqrt{L/\ell}-1)(f(\x_k)-f(\y_{k+1})-\frac{\ell}{2}\Vert \y_{k+1}-\x_k\Vert^2) 
\\
&\hphantom{=}\quad\mbox{}
+ 
\frac{\Vert\nabla f(\y_{k+1})\Vert^2}{2\sqrt{L\ell}}
\quad
\mbox{(by \eref{eq:strcvx})} \\
&\le
f(\x^*)-f(\y_{k+1})-\frac{\ell}{2}\Vert \y_{k+1}-\x^*\Vert^2
+ (\sqrt{L/\ell}-1)(f(\x_k)-f(\y_{k+1})) 
+ 
\frac{\Vert\nabla f(\y_{k+1})\Vert^2}{2\sqrt{L\ell}}
 \\
&=
f(\x^*)-\sqrt{L/\ell}\cdot f(\y_{k+1}) 
-\frac{\ell}{2}\Vert \y_{k+1}-\x^*\Vert^2
+ (\sqrt{L/\ell}-1)f(\x_k) 
+ 
\frac{\Vert\nabla f(\y_{k+1})\Vert^2}{2\sqrt{L\ell}}
\\
&\le 
f(\x^*)-
\sqrt{L/\ell}\cdot
\left(f\left(\y_{k+1}-\nabla f(\y_{k+1})/L\right)+\Vert\nabla  f(\y_{k+1})
\Vert^2/(2L)\right) && \\
&\hphantom{=}\quad\mbox{}
-\frac{\ell}{2}\Vert \y_{k+1}-\x^*\Vert^2
+ (\sqrt{L/\ell}-1)f(\x_k) 
+ 
\frac{\Vert\nabla f(\y_{k+1})\Vert^2}{2\sqrt{L\ell}}
\quad\mbox{(by \eref{eq:cvxdesc})} \\
&=
f(\x^*)-\sqrt{L/\ell}\cdot
f\left(\y_{k+1}-\nabla f(\y_{k+1})/L\right)
-\frac{\ell}{2}\Vert \y_{k+1}-\x^*\Vert^2
+ (\sqrt{L/\ell}-1)f(\x_k) 
 \\
&=
f(\x^*)-\sqrt{L/\ell}\cdot
f(\x_{k+1})
-\frac{\ell}{2}\Vert \y_{k+1}-\x^*\Vert^2
+ (\sqrt{L/\ell}-1)f(\x_k) 
\quad\mbox{(by \eref{eq:xupd2})} \\
&=
(\sqrt{L/\ell}-1)(f(\x_k) -f(\x^*))-\sqrt{L/\ell}\cdot(f(\x_{k+1})-f(\x^*))
-\frac{\ell}{2}\Vert \y_{k+1}-\x^*\Vert^2
\end{align*}

We now can finally analyze the bound on the first term of $C\Psi_{k+1}-\Psi_k$.
We have an explicit formula for \eref{eq:tu1},
we have forced \eref{eq:tu2} to be 0 by choice of $\xi$, and
\eref{eq:tu3ag}
is nonpositive.  The remaining terms are captured in \eref{eq:threeterm}
and \eref{eq:t4bound}, so therefore
\begin{align}
\Vert\xi\w_{k+1}\Vert^2 - \Vert\w_{k}\Vert^2
&\le 
\frac{2}{(\sqrt{L/\ell}-1)\ell}
\notag \\
&\hphantom{\le}\quad\mbox{}
\cdot \left[(\sqrt{L/\ell}-1)(f(\x_k)-f(\x^*))-\sqrt{L/\ell}\cdot
(f(\x_{k+1})-f(\x^*))\right]
\notag \\
&\hphantom{\le}\quad\mbox{}
+\left[\frac{1}{\sqrt{L/\ell}-1} - 
\frac{2}{(\sqrt{L/\ell}-1)\ell}\cdot\frac{\ell}{2}\right]
\Vert\y_{k+1}-\x^*\Vert^2.
\label{eq:xiwdiff}
\end{align}
Observe that the square-bracketed coefficient at the end of 
\eref{eq:xiwdiff} is 0.  Rearranging,
\begin{align*}
\Vert\xi\w_{k+1}\Vert^2 + 
\frac{2}{(\sqrt{L/\ell}-1)\ell}\cdot
\sqrt{L/\ell}\cdot
(f(\x_{k+1})-f(\x^*)) 
&\le 
\Vert\w_{k}\Vert^2\\
&\hphantom{=}\quad\mbox{}
 +  
\frac{2}{(\sqrt{L/\ell}-1)\ell}\\
&\hphantom{=}\quad\quad\mbox{}
\cdot
(\sqrt{L/\ell}-1)(f(\x_k)-f(\x^*)),
\end{align*}
i.e.,
$$\xi^2\left[\Vert\w_{k+1}\Vert^2+\frac{2}{\ell}\cdot
(f(\x_{k+1})-f(\x^*))\right] \le
\Vert\w_{k}\Vert^2 + 
\frac{2}{\ell}\cdot (f(\x_k)-f(\x^*)).$$
Thus, we have established \eref{eq:potdecr}
in the case $C=\xi^2$, which by \eref{eq:xirecurag}
implies
\begin{equation}
C=
1+\frac{1}{\sqrt{L/\ell}-1}.
\label{eq:Cag}
\end{equation}

The case $k=0$ needs separate attention.
First, by taking $\xi=1$, $\rho_0=0$, $\rho_1=\sqrt{L/\ell}-1$,
$\theta_0=\nu_0=0$, we observe 
that  \eref{eq:tu1},
\eref{eq:tu2}, \eref{eq:tu3} and \eref{eq:tu5}
all vanish.
Thus,
\begin{align*}
\Vert\w_1\Vert^2&=
\Vert\w_0\Vert^2+\frac{2(\x^*-\y_1)^T\nabla f(\y_1)}
{\sqrt{L\ell}}
+\frac{\Vert\nabla f(\y_1)\Vert^2}{L\ell} \\
&\le
\Vert\w_0\Vert^2
+\frac{\Vert\nabla f(\y_1)\Vert^2}{L\ell}
\end{align*}
since the dropped term in the last line is nonpositive by convexity.

On the other hand,
\begin{align*}
\frac{2(f(\x_1)-f(\x^*))}{\ell}
&=
\frac{2(f(\y_1-\nabla f(\y_1)/L)-f(\x^*))}{\ell} \\
&\le
\frac{2(f(\y_1)-\Vert\nabla f(\y_1)\Vert^2/(2L)-f(\x^*))}{\ell} 
&\mbox{(by \eref{eq:cvxdesc})}
\\
&=
\frac{2(f(\x_0)-\Vert\nabla f(\y_1)\Vert^2/(2L)-f(\x^*))}{\ell} \\
&=
\frac{2(f(\x_0)-f(\x^*))}{\ell}-\frac{\Vert\nabla f(\y_1)\Vert^2}{L\ell}.
\end{align*}
Adding the two preceding inequalities shows that $\Psi_1\le \Psi_0$.
Thus, in this particular case, \eref{eq:potdecr} does not necessarily
hold for any $C>1$.

\section{Convergence of conjugate gradient}
\label{sec:cvgcg}

We introduce the notation $F_k=2(f(\x_k)-f(\x^*))$.
We define
\begin{equation}
\rho_k=\frac{F_k}{\alpha_k\Vert\r_{k-1}\Vert^2},
\label{eq:rhodeflcg}
\end{equation}
for $k=1,2,\ldots$ and $\rho_0=$.
For $k=1,2,\ldots$, this
$\rho_k$ has the special property that it is the optimizer of
the optimization problem
$\min\{\Vert \x_k+\rho\s_k-\x^*\Vert: \rho\in\R\}$, a property
proved by \cite{hestenesstiefel} (see (6:8)).  This property
is not directly used in the upcoming analysis.

We also require the following result:
\begin{equation}
F_k-F_{k+1}= \alpha_{k+1}\Vert\r_k\Vert^2,
\label{eq:cgdesc}
\end{equation}
which follows from \eref{eq:alphadef},
\eref{eq:xupd}, \eref{eq:prorth}
and the fact that 
\begin{equation}
f(\x_k+\d)=f(\x_k)-\r_k^T\d+\d^TA\d/2
\label{eq:fd}
\end{equation}
for any $\d$.  It is also proven in \cite[(6:1)]{hestenesstiefel}.

We use two other equations from \cite{hestenesstiefel}, the first of
which is (6:5):
\begin{equation}
\Vert \x_k-\x^*\Vert^2-\Vert\x_{k+1}-\x^*\Vert^2 =\frac{(F_k+F_{k+1})\Vert\p_{k+1}\Vert^2}
{\p_{k+1}^TA\p_{k+1}},
\label{eq:hs6:5}
\end{equation}
for $k=0,1,2\ldots$.
Let us assume now that $k\ge 1$; the $k=0$ case is considered below.
The next equation, which holds for $k=1,2,\ldots$, is 
an unnumbered equation of \cite[p.~417, col.~2]{hestenesstiefel}:
\begin{equation}
\Vert \x_k-\x^*\Vert^2 -\Vert\x_k+\rho_k\s_k-\x^*\Vert^2 =
\frac{F_k^2\Vert\p_k\Vert^2}{\Vert\r_{k-1}\Vert^4}.
\label{eq:HSunn}
\end{equation}
If we subtract \eref{eq:hs6:5} and the $k+1$ case of
\eref{eq:HSunn} from the $k$ case of \eref{eq:HSunn}, and
recalling the notation
$\w_k=\x_k+\rho_k\s_k-\x^*$, we obtain
\begin{equation}
\Vert\w_{k+1}\Vert^2 -
\Vert\w_{k}\Vert^2 =
z_1+z_2+z_3+z_4
\end{equation}
where
\begin{align*}
z_1&= \frac{F_k^2\Vert\p_k\Vert^2}{\Vert\r_{k-1}\Vert^4}, \\
z_2&= -\frac{F_{k+1}^2\Vert\p_{k+1}\Vert^2}{\Vert\r_k\Vert^4}, \\
z_3&= -\frac{F_k\Vert\p_{k+1}\Vert^2}{\p_{k+1}^TA\p_{k+1}}, \\
z_4&= -\frac{F_{k+1}\Vert\p_{k+1}\Vert^2}{\p_{k+1}^TA\p_{k+1}}.
\end{align*}
In order to simplify this sum, we
make the following substitutions:
\begin{align*}
F_{k+1}&:=F_{k}-\alpha_{k+1}\Vert\r_k\Vert^2 &\mbox{(by \eref{eq:cgdesc})}, \\
\Vert\r_{k-1}\Vert^2 &:=\Vert\r_k\Vert^2/\beta_{k+1}&\mbox{(by \eref{eq:betadef})}, \\
\Vert\p_k\Vert^2 &:= (\Vert\p_{k+1}\Vert^2-\Vert\r_k\Vert^2)/\beta_{k+1}^2&
\mbox{(by \eref{eq:pupd} and \eref{eq:prorth})}
\end{align*}
to obtain:
\begin{align*}
z_1&= 
\frac{F_k^2(\Vert\p_{k+1}\Vert^2-\Vert\r_k\Vert^2)}
{\Vert\r_k\Vert^4}, \\
z_2&= -\frac{(F_{k}^2-2F_k\alpha_{k+1}\Vert\r_k\Vert^2+\alpha_{k+1}^2\Vert\r_k\Vert^4)\Vert\p_{k+1}\Vert^2}{\Vert\r_k\Vert^4}, \\
z_3&= -\frac{F_k\Vert\p_{k+1}\Vert^2}{\p_{k+1}^TA\p_{k+1}}, \\
z_4&= -\frac{(F_{k}-\alpha_{k+1}\Vert\r_k\Vert^2)\Vert\p_{k+1}\Vert^2}{\p_{k+1}^TA\p_{k+1}}.
\end{align*}
Now let us combine these terms, noting that the first term in $z_1$
cancels the first in $z_2$, and substituting 
$\alpha_{k+1}:=\Vert\r_k\Vert^2/(\p_{k+1}^TA\p_{k+1})$ 
(by \eref{eq:alphadef}) 
to obtain
\begin{align}
\Vert\w_{k+1}\Vert^2 -
\Vert\w_{k}\Vert^2 &=
-\frac{F_k^2}
{\Vert\r_k\Vert^2}
&\mbox{(term from $z_1$)} 
\notag\\
&\hphantom{=}\quad\mbox{}
+
\frac{2F_k\Vert\p_{k+1}\Vert^2}
{\p_{k+1}^TA\p_{k+1}}
-\frac{\Vert\p_{k+1}\Vert^2\cdot\Vert\r_k\Vert^4}
{(\p_{k+1}^TA\p_{k+1})^2}
&\mbox{(terms from $z_2$)}
\notag \\
&\hphantom{=}\quad\mbox{}
-\frac{F_k\Vert\p_{k+1}\Vert^2}{\p_{k+1}^TA\p_{k+1}}
&\mbox{(term from $z_3$)}
\notag \\
&\hphantom{=}\quad\mbox{}
-\frac{F_k\Vert\p_{k+1}\Vert^2}{\p_{k+1}^TA\p_{k+1}}
+\frac{\Vert\r_k\Vert^4\cdot\Vert\p_{k+1}\Vert^2}
{(\p_{k+1}^TA\p_{k+1})^2}
&\mbox{(terms from $z_4$)} 
\notag\\
&=
-\frac{F_k^2}
{\Vert\r_k\Vert^2}.
\label{eq:wdiffcg}
\end{align}
It is also possible to obtain \eref{eq:wdiffcg} from
\eref{eq:tu1}--\eref{eq:tu6} with the choice $\xi=1$.

Another helpful inequality is
\begin{align}
\Vert \w_k\Vert^2 &\le
\Vert \x_k-\x^*\Vert^2 
& \mbox{(by \eref{eq:HSunn})}
\notag\\
&\le (\x_k-\x^*)^T(A/\ell)(\x_k-\x^*) 
& \mbox{(since $\lambda_{\min}(A/\ell)\ge 1$)} 
\notag \\
& = F_k/\ell. \label{eq:wbd1}
\end{align}

Now we establish \eref{eq:potdecr}:
\begin{align*}
C\Psi_{k+1}-\Psi_k &=
C(\Vert\w_{k+1}\Vert^2+F_{k+1}/\ell) 
- (\Vert \w_{k}\Vert^2+F_k/\ell) \\
&=(C-1)(\Vert\w_{k+1}\Vert^2+F_{k+1}/\ell)
+\Vert\w_{k+1}\Vert^2-\Vert\w_k\Vert^2 
\\
&\hphantom{=}\quad\mbox{}
+ (F_{k+1}-F_k)/\ell \\
&=
(C-1)(\Vert\w_{k+1}\Vert^2+F_{k+1}/\ell)
-F_k^2/\Vert\r_k\Vert^2
\\
&\hphantom{=}\quad\mbox{}
-\alpha_{k+1}\Vert\r_k\Vert^2/\ell 
&\mbox{(by \eref{eq:cgdesc} and \eref{eq:wdiffcg})}
\\
&\le
(C-1)(\Vert\w_{k+1}\Vert^2+F_{k+1}/\ell)
-F_k^2/\Vert\r_k\Vert^2
\\
&\hphantom{=}\quad\mbox{}
-\Vert\r_k\Vert^2/(L\ell ) 
&\mbox{(by \eref{eq:alphaRQ})}
\\
&\le 
(C-1)(\Vert\w_{k+1}\Vert^2+F_{k+1}/\ell)
-2F_k/\sqrt{L\ell} 
&\mbox{(since $x^2+y^2\ge 2xy$)}
\\
&\le
2(C-1)F_{k+1}/\ell - 2F_k/\sqrt{L\ell} 
&\mbox{(by \eref{eq:wbd1})}
\\
&\le
2(C-1)F_k/\ell - 2F_k/\sqrt{L\ell} \\
&\le 0
\end{align*}
provided that $2(C-1)/\ell\le 2/\sqrt{L\ell}$, i.e.,
\begin{equation}
C\le 1+\sqrt{\frac{\ell}{L}}.
\label{eq:Ccg}
\end{equation}
Thus, we take $C$ equal to the right-hand side  of the preceding
inequality to establish \eref{eq:potdecr}.

Again, the $k=0$ case needs special attention.  For this case,
as with accelerated gradient, we settle for the weaker inequality
that $\Psi_1\le \Psi_0$.  This inequality holds for each of the
two terms separately:
\begin{align*}
\Vert \w_1\Vert^2 &\le \Vert \x_1-\x^*\Vert^2
&\mbox{(by \eref{eq:HSunn})} \\
&\le \Vert \x_0-\x^*\Vert^2
&\mbox{(by \eref{eq:hs6:5})} \\
& = \Vert\w_0\Vert^2.
\end{align*}
Also,  $F_1\le F_0$ by \eref{eq:cgdesc}.  

\section{Discussion}
\label{sec:discussion}

The main point of this work is to prove Theorem~\ref{mainthm}
using the same convergence bound for accelerated gradient and conjugate
gradient.  The result in this paper was originally motivated by
our consideration of nonlinear conjugate gradient.

The traditional extensions of linear CG to nonlinear CG
for the general
case of unconstrained optimization, e.g.,
the algorithms of  Fletcher and Reeves \cite{fletcherreeves} 
and Polak and Ribi\`ere \cite{polakribiere} (see
Nocedal and Wright \cite{NocedalWright} for an overview of these
algorithms)
are not optimal
for minimizing strongly convex functions.  
Unlike accelerated gradient, there is no global complexity bound known
for any nonlinear CG method even in the case of strongly convex functions.
Indeed, Nemirovsky
and Yudin argue that traditional nonlinear CG can perform even worse than
steepest descent.  

A partial unification the analyses of linear CG and accelerated gradient such
as ours could point the way to development of a new nonlinear CG
method.  The new method would have two desirable properties: (1) it reduces to
linear CG in the case of a quadratic function, and (2) it
maintains the global convergence bound of accelerated gradient.
Furthermore, such an algorithm would ideally be able to adapt
between steps of the two algorithms even within the same problem.
A preliminary proposal for a nonlinear CG like this
was made in the PhD thesis
of the first author \cite{Karimi:thesis}, and will be
the subject of ongoing work.

A second practical use of the unified analysis is the
consideration of algorithms for minimizing a quadratic
function using a modification of
linear CG, such that the modification 
changes it into
a nonlinear iteration.  For example, several authors \cite{SimonciniSzyld,
VdeSleijpen} have considered the use of conjugate gradient methods
in the case of noisy matrix-vector multiplication.  Other
authors, e.g., \cite{Knyazev} have considered the possibility
of changing the preconditioner from one iterate to the next.
In both of these cases, optimality with respect to the Krylov
space is no longer assured.  However, it is  possible
that the bound in Theorem~\ref{mainthm} may still hold.
Since the analysis establishes Theorem~\ref{mainthm} without
relying on Krylov optimality, it may enable new analyses of
such `perturbed' conjugate gradient methods.  This matter
is also left for future work.

Another application of the LCG bound developed herein is to computational
scientists developing new linear conjugate gradient methods (e.g.,
new preconditioners or new ways to compute matrix-vector products).
Our bound not only directly shows the convergence rate
claimed in Theorem~\ref{mainthm}, but more strongly it
shows that the potential decreases by at least a fixed
constant factor on each iteration. 
In a test run of any proposed new algorithm, it is possible to
measure the potential developed herein and monitor its
steady decrease.  Any failure to exhibit the prescribed decrease would be an unambiguous
indication that the method is failing due to some source of inexactness
(e.g., roundoff error).  In contrast, better known measures of LCG
convergence can stagnate for many consecutive iterations, making it
difficult to detect the impact of inexactness.  We remark that
in order to use our potential in this manner, 
it is of course necessary
to know $\ell$, $L$ and the exact solution to the linear system at the outset, which
is often the case in testing a new algorithm but obviously not in
its practical use.

As for the theoretical content of this paper, it would be useful
to simplify our analysis, which appears to complicated, and also
to further unify the treatments of the two algorithms.
Another useful development would be a potential  that
involves the term $\Vert\nabla f(\x_k)\Vert^2$.  This is because,
in practice, an algorithm does not have access to either 
$f(\x_k)-f(\x^*)$ or $\Vert\x_k-\x^*\Vert$, so the potential proposed
here could not be evaluated by an algorithm to measure progress.

\bibliography{optimization}
\bibliographystyle{plain}
\end{document}